\documentclass[final,1p,times]{elsarticle}

\usepackage{amssymb}
 \usepackage{amsthm}
\usepackage{amscd}
\usepackage{bm}
\usepackage{amsmath}
\usepackage{amsfonts}
\usepackage{amssymb}
\usepackage{graphicx}
\newtheorem{theorem}{Theorem}
\newtheorem{proposition}[theorem]{Proposition}
\newtheorem{lemma}[theorem]{Lemma}
\newtheorem{corollary}[theorem]{Corollary}
\usepackage{mathrsfs}
\usepackage{titletoc}

\newcommand{\ra}{\rightarrow}
\newcommand{\p}{\partial}
\newcommand{\f}{\frac}

\newcommand{\be}{\begin{equation}}
\renewcommand{\ra}{\rightarrow}
\newcommand{\ee}{\end{equation}}
\newcommand{\bea}{\begin{eqnarray}}
\newcommand{\eea}{\end{eqnarray}}
\newcommand{\bna}{\begin{eqnarray*}}
\newcommand{\ena}{\end{eqnarray*}}

\renewcommand{\le}{\left}
\newcommand{\ri}{\right}

\journal{***}

\begin{document}

\begin{frontmatter}

\title{A gradient flow for the prescribed Gaussian curvature problem on a closed Riemann surface with conical singularity}

 \author{Yunyan Yang}
 \ead{yunyanyang@ruc.edu.cn}

\address{ Department of Mathematics,
Renmin University of China, Beijing 100872, P. R. China}

\begin{abstract}
 In this note, we prove that the abstract gradient flow introduced by Baird-Fardoun-Regbaoui \cite{BFR}
 is well-posed on a closed Riemann surface with conical singularity. Long time existence and convergence
 of the flow are proved under certain assumptions. As an application, the prescribed Gaussian curvature problem is 
 solved when the singular Euler characteristic of the conical surface is non-positive.  
\end{abstract}


\end{frontmatter}

\titlecontents{section}[0mm]
                       {\vspace{.2\baselineskip}}
                       {\thecontentslabel~\hspace{.5em}}
                        {}
                        {\dotfill\contentspage[{\makebox[0pt][r]{\thecontentspage}}]}
\titlecontents{subsection}[3mm]
                       {\vspace{.2\baselineskip}}
                       {\thecontentslabel~\hspace{.5em}}
                        {}
                       {\dotfill\contentspage[{\makebox[0pt][r]{\thecontentspage}}]}

\setcounter{tocdepth}{2}


\section{Introduction}

Let $\Sigma$ be a closed Riemann surface, $g$ be a smooth metric and $\kappa$ be its Gaussian curvature. If $\tilde{g}=e^{2u}g$
for some smooth function $u$, then the Gaussian curvature of $\tilde{g}$ satisfies
 $\tilde{\kappa}=e^{-2u}(\Delta_gu+\kappa)$,
 where $\Delta_g$ is the Laplace-Beltrami operator. For a given function $K:\Sigma\ra\mathbb{R}$, can one find a
 metric $\tilde{g}=e^{2u}g$ having $K$ as its Gaussian curvature? This problem is equivalent to the solvability of the equation
 \be\label{K-Z}\Delta_gu+\kappa-Ke^{2u}=0.\ee
 Integration by parts and the Gauss-Bonnet formula imply that necessarily $K$ must have the same sign as the
 topological Euler characteristic $\chi(\Sigma)$ somewhere and in the case $\chi(\Sigma)=0$, either $K$ is identically zero or changes sign.
 It is natural to ask if this condition is also sufficient to guarantee a solution.

 In the case $\chi(\Sigma)<0$, via the method of upper and lower solutions, it was shown by Kazdan-Warner \cite{K-W} that
 if $K\leq 0$ and $K\not\equiv 0$, then (\ref{K-Z}) has a solution. Suppose that $K\leq \sup_\Sigma K=0$,
$K\not\equiv 0$, and $\lambda\in\mathbb{R}$. Using a variational method, Ding-Liu
\cite{Ding-Liu} proved the following: Replacing $K$ by $K+\lambda$ in (\ref{K-Z}), one finds some constant $\lambda^\ast>0$ such that  if $0<\lambda<\lambda^\ast$, then (\ref{K-Z}) has at least two different solutions;
if $\lambda=\lambda^\ast$, then (\ref{K-Z}) has at least one solution; while if $\lambda>\lambda^\ast$, then (\ref{K-Z}) has no solution.
 In the case $\chi(\Sigma)=0$, the problem was completely solved. It was proved by Berger \cite{Berger} that
 if $K\equiv 0$ or $K$ changes sign and $\int_\Sigma Ke^{2v}dv_g<0$,
 where $v$ is a solution of $\Delta_gv=-\kappa$, then (\ref{K-Z}) has a solution. Later Kazdan-Warner \cite{K-W} pointed out that
 the above assumptions on $K$ is also necessary.
 If $\chi(\Sigma)>0$, $\Sigma$ is either the projective space $\mathbb{RP}^2$
or the $2$-sphere $\mathbb{S}^2$. In the case of $\mathbb{RP}^2$,
it was shown by Moser \cite{Moser} that (\ref{K-Z}) has a solution provided that
$\sup_\Sigma K>0$ and $K(p)=K(-p)$ for all $p\in \mathbb{S}^2$. While the problem on $\mathbb{S}^2$ is much more
complicated and known as the Nirenberg problem. Moser's result was extended by Chang-Yang \cite{Chang-Yang} to reflected symmetric function
$K$ under further assumptions. For rotationally symmetric function $K$, sufficient condition was given by Chen-Li \cite{Chen-Li} and
Xu-Yang \cite{Xu-Yang}. Concerning more general functions $K$, we refer the reader to \cite{Chang-Yang2,CGY,Chen-Ding}.

Also various flows have been employed to attack the problem.
In \cite{Hamilton}, The Ricci flow was introduced by Hamilton to find a solution of (\ref{K-Z}), where $K$ is a constant.
His result was later completed by Chow
\cite{Chow}. The Calabi flow was investigated by Bartz-Struwe-Ye \cite{B-S-Y} and Struwe \cite{St}.
While in \cite{Str}, Struwe used the Gaussian curvature flow to reprove Chang-Yang's results \cite{Chang-Yang}.
For further developments of this flow, we refer the reader to
Brendle \cite{Brendle1,Brendle2}, Ho \cite{Ho} and Zhang \cite{Zhang}.
Assuming that the initial metric $g$ has constant Gaussian curvature $\kappa$. Baird-Fardoun-Regbaoui \cite{BFR}
proposed an abstract gradient flow,
through which $g(t)$ converges to a metric having the prescribed Gaussian curvature. This method
solved (\ref{K-Z}) perfectly in the case $\chi(\Sigma)\leq 0$ and partially in the case $\chi(\Sigma)>0$.

The same problem can be proposed on conical surfaces. We begin with basic definitions.
Let $\Sigma$ be a closed Riemann surface as before. A metric $g$ is said to be
a conformal metric having conical singularity of order $\beta>-1$ at $p\in\Sigma$, if in a local holomorphic
coordinate with $z(p)=0$, there exists some function $u$ which is continuous and $C^2$ away from zero such that
$$g=e^{2u}|z|^{2\beta}|dz|^2.$$ If $g$ has conical singularities of order $\beta_i>-1$ at $p_i\in\Sigma$, $i=1,\cdots,\ell$, we say that
$g$ represents a divisor $\bm\beta=\sum_{i=1}^\ell \beta_ip_i$. Then the pair $(\Sigma,\bm\beta)$
is called a conical surface, and the corresponding singular Euler characteristic
is written as
\be\label{singular-Euler}\chi(\Sigma,\bm\beta)=\chi(\Sigma)+\sum_{i=1}^\ell\beta_i,\ee
where $\chi(\Sigma)$ is the topological Euler characteristic.

If $\chi(\Sigma,\bm\beta)$ is nonpositive, the problem can be solved in the variational framework as the case of smooth metrics. Precisely,
it was shown by Troyanov \cite{Troyanov} that if $\chi(\Sigma,\bm\beta)<0$, then any smooth negative function is the Gaussian curvature
of a unique conformal metric $\tilde{g}$ representing $\bm\beta$. Recently this result has been improved by Zhu and
the author \cite{Yang-Zhu} by using the variational method of Ding-Liu \cite{Ding-Liu} and Borer-Galimberti-Struwe \cite{B-G-Stru}.
In particular, if we assume $\chi(\Sigma,\bm\beta)<0$, the background metric $g$ has the Gaussian curvature $\kappa\equiv -1$,
and $K$ is a smooth function satisfying $\sup_\Sigma K=0$ and $K\not\equiv 0$, then there exists
a unique function
$$u\in \mathscr{C}=C^2(\Sigma\setminus\{p_1,\cdots,p_\ell\})\cap C^0(\Sigma)\cap W^{1,2}(\Sigma,g)$$
such that the metric $e^{2u}g$ has the Gaussian curvature $K$; moreover, there exists
some constant $\lambda^\ast>0$ such that when $0<\lambda<\lambda^\ast$, there exist at least two different functions
$u_1,u_2\in \mathscr{C}$
such that $e^{2u_1}g$ and $e^{2u_2}g$ have the same Gaussian curvature
$K+\lambda$; when $\lambda=\lambda^\ast$, there exists at least one function $u\in\mathscr{C}$ such that $e^{2u}g$
has the Gaussian curvature $K+\lambda^\ast$; when $\lambda>\lambda^\ast$, there is no function $u\in W^{1,2}(\Sigma,g)$
such that $e^{2u}g$ has the Gaussian curvature $K+\lambda$.
The problem was completely solved by Troyanov \cite{Troyanov} in the case $\chi(\Sigma,\bm\beta)=0$. Namely,
there exists a flat metric $g$ representing
$\bm\beta$; moreover, a smooth function $K$ is the Gaussian curvature of a metric $\tilde{g}$ conformal to $g$
if and only if $K$ changes sign and $\int_\Sigma Kdv_{g}<0$.
If $\chi(\Sigma,\bm\beta)>0$, then the problem becomes very subtle. There is much interesting work concerning this situation,
see for examples Troyanov \cite{Troyanov}, McOwen \cite{McOwen}, Chen-Li \cite{Chen-Li1,Chen-Li2,chen-Li-13}, Luo-Tian \cite{Luo-Tian},
Mondello-Panov \cite{Mondello-Panov}, Bartolucci \cite{Bart}, Bartolucci-De Marchis-Malchiodi \cite{B-D-M}, Fang-Lai
\cite{F-L} and a very nice survey of Lai \cite{Lai}.

Again the Ricci flow is an elegant way to solve the problem on conical surfaces. Yin \cite{Yin1,Yin2,Yin3} established
a basic theory in this regards, and proved long time existence and convergence of the flow when $\chi(\Sigma,\bm\beta)\leq 0$.
The convergence in the case $\chi(\Sigma,\bm\beta)>0$ was studied by Phong-Song-Sturm-Wang \cite{PSSW1,PSSW2}.
Another approach for the Ricci flow was proposed by Mazzeo-Rubinstein-Sesum \cite{MRS}.

Our aim is to establish the gradient flow of Baird-Fardoun-Regbaoui \cite{BFR} on conical surfaces.
Assuming the background metric has a constant Gaussian curvature, we prove the long time existence of the flow.
Moreover, when  $\chi(\Sigma,\bm\beta)\leq 0$, we obtain the convergence of the flow under additional assumptions.
For the proof of our results, we follow the lines of
Baird-Fardoun-Regbaoui \cite{BFR}.
Here the key point is the following observation:  the functionals involved
are still analytic if the background metric  has conical singularity.

The remaining part of this note is organized as follows: In Section 2, we construct functional framework
and give main results of this note; In Section 3, we prove the analyticity of functionals $\mathcal{J}$ and $\mathcal{L}$,
and calculate their gradients; In Section 4, we show the long time existence of the gradient flow; In Section 5, a sufficient condition
for convergence of the flow will be discussed; In Section 6, we prove that when $\chi(\Sigma,\bm\beta)\leq 0$, the flow converges to the desired
solution of the problem.

\section{Notations and main results}
Let $\Sigma$ be a closed Riemann surface, ${\bm\beta}=\sum_{i=1}^\ell\beta_ip_i$
be a divisor, $\beta_i>-1$ for all $i$, and $g$ be a conformal
metric representing $\bm\beta$.
Let $\kappa: \Sigma\setminus {\rm supp}\,
\bm\beta\ra \mathbb{R}$ be the Gaussian curvature of $g$, where
${\rm supp}\,\bm\beta=\{p_1,\cdots,p_\ell\}$. From now on, we assume $\kappa$ is a constant.
Then the Gauss-Bonnet formula (see for example \cite{Troyanov}) reads
$$\label{Gauss-Bonnet}\kappa {\rm Vol}_g(\Sigma)=\int_\Sigma \kappa dv_g=2\pi\chi(\Sigma,\bm\beta),$$
where $\chi(\Sigma,\bm\beta)$ is defined as in (\ref{singular-Euler}), and $dv_g$ denotes the volume element with respect to the conical metric $g$.
Clearly there exists a smooth metric $g_0$ such that
$$\label{g0}g=\rho g_0,$$ where
$\rho>0$ on $\Sigma$, $\rho\in C^2_{\rm loc}(\Sigma\setminus{\rm supp}\,\bm\beta)$,
and $\rho\in L^{r}(\Sigma)$ for some $r>1$.
 Let $W^{1,2}(\Sigma,g)$ be the completion of $C^\infty(\Sigma)$ under the norm
$$\|u\|_{W^{1,2}(\Sigma,g)}=\le(\int_\Sigma(|\nabla_gu|^2+u^2)dv_g\ri)^{1/2},$$
where $\nabla_g$ denotes the gradient operator with respect to the metric $g$.
It was observed by Troyanov \cite{Troyanov} that $W^{1,2}(\Sigma,g)=W^{1,2}(\Sigma,g_0)$.
In particular,
$W^{1,2}(\Sigma,g)$ is a Hilbert space, which is hereafter denoted by $\mathscr{H}$, with an inner product
$$\langle u,w\rangle_{\mathscr{H}}=\int_\Sigma (\nabla_gu\nabla_gw+uw)dv_g.$$
 Moreover, by the Sobolev embedding theorem for
smooth Riemann surface $(\Sigma,g_0)$ and the H\"older inequality, one has
$$\label{Sobolev}W^{1,2}(\Sigma,g)\hookrightarrow L^p(\Sigma,g),\quad \forall p>1.$$
Let $\bar{g}=e^{2u}g$ be another conical metric representing $\bm\beta$ and $K:\Sigma\setminus{\rm supp}\bm\beta
\ra\mathbb{R}$ be the Gaussian curvature of $\bar{g}$.
Then $K$ satisfies point-wisely on $\Sigma\setminus{\rm supp}\bm\beta$,
 \be\label{pres-equ}K=e^{-2u}({\kappa}+\Delta_gu),\ee
 where $\Delta_g$ denotes the Laplacce-Beltrami
operator with respect to the metric $g$. Obviously, if $u$ is a distributional solution of the equation
\be\label{Kz-eqn}\Delta_gu+\kappa-Ke^{2u}=0,\ee
then $u$ satisfies (\ref{pres-equ}).

 Let us define two functionals
$\mathcal{J}:\mathscr{H}\ra\mathbb{R}$, $\mathcal{L}:\mathscr{H}\ra\mathbb{R}$ by
\bea\label{funct}&&\mathcal{J}(u)=\f{1}{2}\int_\Sigma|\nabla_gu|^2dv_g+\kappa\int_\Sigma udv_g,\\
\label{funct-L}&&\mathcal{L}(u)=\f{1}{2}\int_\Sigma Ke^{2u}dv_g,\eea
and a set of functions by
\be\label{constraint}\bm{\mathscr{S}}=\le\{u\in \mathscr{H}: \mathcal{L}(u)=\kappa{\rm Vol}_g(\Sigma)=2\pi\chi(\Sigma,\bm\beta)\ri\}.\ee
The gradients of $\mathcal{J}$ and $\mathcal{L}$,
$\nabla\mathcal{J}:\mathscr{H}\ra\mathscr{H}$ and $\nabla\mathcal{L}:\mathscr{H}\ra\mathscr{H}$ are defined by
\bea\label{nab-j}&&\langle\nabla\mathcal{J}(u),w\rangle_{\mathscr{H}}=d\mathcal{J}(u)(w)=\le.\f{d}{dt}\ri|_{t=0}\mathcal{J}(u+tw),\\
\label{nab-l}&&\langle\nabla\mathcal{L}(u),w\rangle_{\mathscr{H}}=d\mathcal{L}(u)(w)=\le.\f{d}{dt}\ri|_{t=0}\mathcal{L}(u+tw)\eea
respectively, where $u$ and $w$ are functions taken from $\mathscr{H}$.
Hereafter we assume $K\not\equiv 0$. It follows that $\nabla \mathcal{L}(u)\not=0$ for all $u\in{\mathscr{S}}$.
Thus ${\mathscr{S}}$ is a smooth hypersurface in
$\mathscr{H}$. A unit normal on $\mathscr{S}$ is
$$\label{normal}\mathcal{N}(u)=\f{\nabla \mathcal{L}(u)}{\|\nabla \mathcal{L}(u)\|_{\mathscr{H}}}$$
for any $u\in\mathscr{S}$, where $\|\cdot\|_{\mathscr{H}}=\langle\cdot,\cdot\rangle_{\mathscr{H}}$.
This allows us to consider the gradient of $\mathcal{J}$ with respect to
the hypersurface $\mathscr{S}$, which is defined by
\be\label{nab-C}\nabla^{{\mathscr{S}}}\mathcal{J}(u)=\nabla \mathcal{J}(u)-\langle\nabla \mathcal{J}(u),\mathcal{N}(u)\rangle_\mathscr{H}
\mathcal{N}(u).\ee
The gradient flow of $\mathcal{J}$ with respect to the hypersurface $\mathscr{S}$ can be written as
\be\label{flow}\le\{
\begin{array}{lll}
\p_tu=-\nabla^{{\mathscr{S}}}\mathcal{J}(u)\\[1.2ex]
u(0)=u_0\in{\mathscr{S}}.
\end{array}\ri.\ee
If the flow exists for all time and converges at infinity, then the limit function $u_\infty$ gives
a distributional solution of (\ref{Kz-eqn}). Our first result is an analog of (\cite{BFR}, Theorem 1), namely
\begin{theorem}\label{Theorem 1}
Let $\Sigma$ be a closed Riemann surface, $\bm\beta=\sum_{i=1}^\ell \beta_ip_i$ be a divisor with $\beta_i>-1$, $i=1,\cdots,\ell$,
 and $g$ be a metric representing $\bm\beta$. Let $\mathcal{J}$, $\mathcal{L}$ and $\mathscr{S}$ be defined by (\ref{funct}),
  (\ref{funct-L}) and (\ref{constraint}) respectively. Suppose that the Gaussian curvature of $g$ is a constant $\kappa$,
  and that $K\in C^0(\Sigma)$ satisfies the condition
\be\label{compati}\le\{\begin{array}{lll}
\int_\Sigma Kdv_g<0&{\rm when} &\chi(\Sigma,\bm\beta)<0\\[1.2ex]
\int_\Sigma Kdv_g<0,\,\,\sup_{\Sigma}K>0&{\rm when}&\chi(\Sigma,\bm\beta)=0\\[1.2ex]
\sup_{\Sigma}K>0&{\rm when}&\chi(\Sigma,\bm\beta)>0.
\end{array}
\ri.\ee
Then for any $u_0\in\mathscr{S}$,
there exists a unique global solution $u\in C^\infty([0,\infty), \mathscr{H})$
of the gradient flow (\ref{flow}), satisfying $u(t)\in \mathscr{S}$ for all $t\geq 0$. Moreover the energy identity
\be\label{energy-id}\int_0^t\|\p_su(s)\|^2ds+\mathcal{J}(u(t))=\mathcal{J}(u_0).\ee
holds for all $t>0$.
\end{theorem}
If $\chi(\Sigma,\bm\beta)\leq 0$, then we have the convergence of the flow, an analog of
(\cite{BFR}, Theorem 2).

\begin{theorem}\label{theorem2}
Let $u_0\in\mathscr{S}$ and $u:[0,\infty)\ra\mathscr{H}$ be given as in Theorem \ref{Theorem 1}. In the case $\chi(\Sigma,\bm\beta)=0$,
there exists a $u_\infty\in W^{2,r}(\Sigma,g)\cap C^\alpha(\Sigma)$ for some $r>1$ and $0<\alpha<1$ such that $u(t)$ converges to
$u_\infty$ in $\mathscr{H}$ as $t\ra\infty$, moreover $u_\infty+\tau$ is a distributional solution of (\ref{Kz-eqn}) for some constant $\tau$; In the case
$\chi(\Sigma,\bm\beta)<0$, there exists a positive constant $\epsilon_0$ depending only on $K^-(x)=\max\{-K(x),0\}$ and the conical
metric $g$ such that if $u_0$ satisfies
\be\label{small}e^{\gamma \|u_0\|_{\mathscr{H}}^2}\sup_{x\in\Sigma}K(x)\leq \epsilon_0,\ee
where $\gamma>1$ is a constant depending only on $g$, then $u(t)$ converges in $\mathscr{H}$ to a distributional solution
$u_\infty$ of (\ref{Kz-eqn}) as $t\ra\infty$.
\end{theorem}

We remark that if $K(x)\leq 0$, then the hypothesis (\ref{small}) is obviously satisfied by all $u_0\in\mathscr{H}$.
Finally, as an interesting application of Theorem \ref{theorem2}, we have the following:

\begin{corollary}\label{Cor1}
Suppose $K\in C^0(\Sigma)$ and $\int_\Sigma Kdv_g<0$.
If in addition $\sup_{x\in\Sigma}K(x)>0$ in the case $\chi(\Sigma,\bm\beta)=0$, or $\sup_{x\in\Sigma}\max\{K(x),0\}$ is sufficiently small
in the case $\chi(\Sigma,\bm\beta)<0$,
then  there exists a conformal metric $\tilde{g}$ representing $\bm\beta$ and having $K$ as its Gaussian curvature.
\end{corollary}

\section{Preliminaries}
In this section, we first show the analyticity of the functionals $\mathcal{J}$ and $\mathcal{L}$, and then calculate their gradients.
\begin{lemma}\label{analytic}
The functionals $\mathcal{J}:\mathscr{H}\ra\mathbb{R}$ and $\mathcal{L}:\mathscr{H}\ra\mathbb{R}$ are analytic.
\end{lemma}

\noindent{\it Proof.} Let $u, h\in\mathscr{H}$ be fixed. Clearly $\mathcal{J}$ has the following Taylor expansion (see for example Chang
\cite{Chang}, Theorem 1.4 of Chapter 1)
\be\label{Taylor}\mathcal{J}(u+h)=\sum_{k=0}^n\f{\mathcal{J}^{(k)}(u)h^{(k)}}{k!}+\mathcal{R}_n(u,h)h^{(n)},\ee
where $\mathcal{J}^{(0)}(u)=\mathcal{J}(u)$, $h^{(k)}$ stands for $(\underbrace{h,\cdots,h}_{k})$, $k=0,1,2,\cdots$, and $\mathcal{R}_n(u,h)$ satisfies
\be\label{err-J}\mathcal{R}_n(u,h)=\int_0^1\f{(1-t)^{n-1}}{(n-1)!}\le\{\mathcal{J}^{(n)}(u+th)-\mathcal{J}^{(n)}(u)\ri\}dt.\ee
One easily computes when $n\geq 3$,
\bna &&\mathcal{J}^{(n)}(u)h^{(n)}=\le.\f{\p^n}{\p t_1\cdot\p t_n}\mathcal{J}(u+t_1h+\cdots+t_nh)\ri|_{t_1=\cdots=t_n=0}=0,\\
&&\mathcal{J}^{(n)}(u+th)h^{(n)}=\le.\f{\p^n}{\p t_1\cdot\p t_n}\mathcal{J}(u+th+t_1h+\cdots+t_nh)\ri|_{t_1=\cdots=t_n=0}=0.\ena
Hence we have
\be\label{error}\lim_{n\ra\infty}\mathcal{R}_n(u,h)h^{(n)}=0.\ee
Combining (\ref{Taylor}) and (\ref{error}), we conclude that $\mathcal{J}:\mathscr{H}\ra\mathbb{R}$ is analytic.

Similar to (\ref{Taylor}), we have
\be\label{Taylor-L}\mathcal{L}(u+h)=\sum_{k=0}^n\f{\mathcal{L}^{(k)}(u)h^{(k)}}{k!}+\mathcal{R}_n^\mathcal{L}(u,h)h^{(n)},\ee
where $\mathcal{R}_n^\mathcal{L}(u,h)h^{(n)}$ is an analog of (\ref{err-J}) with $\mathcal{J}$ replaced by $\mathcal{L}$.
 In view of (\ref{funct-L}), we have for all $n\in\mathbb{N}$, $t\in[0,1]$,
\bna&&\mathcal{L}^{(n)}(u)h^{(n)}=\int_\Sigma Ke^{2u}h^ndv_g,\\
&& \mathcal{L}^{(n)}(u+th)h^{(n)}=\int_\Sigma Ke^{2(u+th)}h^ndv_g.\ena
Clearly there holds for all $t\in[0,1]$,
\bea\le|\le(\mathcal{L}^{(n)}(u+th)-\mathcal{L}^{(n)}(u)\ri)h^{(n)}\ri|&\leq& \int_\Sigma |K|e^{2(|u|+|h|)}|h|^ndv_g
\nonumber\\\label{er-l}&\leq&
\le(\int_\Sigma K^2e^{4(|u|+|h|)}dv_g\ri)^{1/2}\le(\int_\Sigma h^{2n}dv_g\ri)^{1/2}.\eea
 It follows that
\bea\nonumber|\mathcal{R}_n^\mathcal{L}(u,h)h^{(n)}|&\leq& \le(\int_\Sigma K^2e^{4(|u|+|h|)}dv_g\ri)^{1/2}\le(\int_\Sigma h^{2n}dv_g\ri)^{1/2}
\f{1}{n!}\\\nonumber&=&\f{1}{\sqrt{n!}}\le(\int_\Sigma K^2e^{4(|u|+|h|)}dv_g\ri)^{1/2}\le(\int_\Sigma \f{h^{2n}}{n!}dv_g\ri)^{1/2} \\
&\leq&\f{1}{\sqrt{n!}}\le(\int_\Sigma K^2e^{4(|u|+|h|)}dv_g\ri)^{1/2}\le(\int_\Sigma e^{h^2}dv_g\ri)^{1/2}\label{R-est}
\eea
Since $u$ and $h$ are fixed functions in $\mathscr{H}$, by a singular Trudinger-Moser inequality (\cite{Troyanov}, Theorem 6), both $e^{2(|u|+|h|)}$ and $e^{u^2}$ belong to $L^p(\Sigma,g)$
for any $p>1$. Note also $K\in C^0(\Sigma)$. Then it follows from (\ref{R-est}) that
$$\lim_{n\ra\infty}\mathcal{R}_n^\mathcal{L}(u,h)h^{(n)}=0.$$
This together with (\ref{Taylor-L}) implies that $\mathcal{L}: \mathscr{H}\ra\mathbb{R}$ is analytic. $\hfill\Box$\\

Let $I$ be an identity operator. We now define a map
$(\Delta_g+I)^{-1}: L^2(\Sigma,g)\ra \mathscr{H}$ in the following way.
For any $f\in L^2(\Sigma,g)$, we say $u=(\Delta_g+I)^{-1}f\in\mathscr{H}$
provided that $(\Delta_g+I)u=f$. Though in our setting, the metric $g$ has conical singularity, the existence and uniqueness of $u$
follows from the Lax-Milgram theorem. Thus the map $(\Delta_g+I)^{-1}$ is well defined.
Moreover $(\Delta_g+I)^{-1}$ is a linear map, which follows from the linearity of $\Delta_g+I$.
Now we have

\begin{lemma}\label{grds}
The gradients of $\mathcal{J}$ and $\mathcal{L}$  at $u\in\mathscr{H}$ are calculated by
\bea\label{gradient-J}&&\nabla \mathcal{J}(u)=u-(\Delta_g+I)^{-1}(u-\kappa),\\[1.2ex]
\label{gradient-2}&&\nabla \mathcal{L}(u)=(\Delta_g+I)^{-1}(Ke^{2u}).\eea
\end{lemma}

\noindent{\it Proof.} On one hand, integration by parts gives
\be\label{onehand}\langle\nabla\mathcal{J}(u),w\rangle_{\mathscr{H}}=\int_\Sigma\le(\nabla_g\nabla\mathcal{J}(u)\nabla_gw+
\nabla\mathcal{J}(u)w\ri)dv_g=\int_\Sigma (\Delta_g+I)\nabla\mathcal{J}(u)wdv_g.\ee
On the other hand,
\be\label{otherhand}d\mathcal{J}(u)(w)=\le.\f{d}{dt}\ri|_{t=0}\mathcal{J}(u+tw)=\int_\Sigma(\nabla_gu\nabla_gw+\kappa w)dv_g=\int_\Sigma(\Delta_gu+\kappa)wdv_g.\ee
Combining (\ref{nab-j}), (\ref{onehand}) and (\ref{otherhand}), we have
$$(\Delta_g+I)\nabla\mathcal{J}(u)=\Delta_gu+\kappa=(\Delta_g+I)u-(u-\kappa),$$
 which leads to
 $$(\Delta_g+I)(\nabla\mathcal{J}(u)-u)=-(u-\kappa).$$
 Then (\ref{gradient-J}) follows immediately.

 To calculate $\nabla\mathcal{L}(u)$, we firstly have an analog of (\ref{onehand}),
$$\langle\nabla\mathcal{L}(u),w\rangle_{\mathscr{H}}=\int_\Sigma (\Delta_g+I)\nabla\mathcal{L}(u)wdv_g.$$
Secondly we have
$$d\mathcal{L}(u)(w)=\le.\f{d}{dt}\ri|_{t=0}\mathcal{L}(u+tw)=\int_\Sigma Ke^{2u}wdv_g.$$
Finally, in view of (\ref{nab-l}), we obtain (\ref{gradient-2}).
$\hfill\Box$

\section{Long time existence and energy identity}
In this section, we prove Theorem \ref{Theorem 1} by following the lines of Baird-Fardoun-Regbaoui \cite{BFR}.\\

\noindent{\it Proof of Theorem \ref{Theorem 1}.} By (\ref{gradient-2}), we have $\nabla\mathcal{L}(u)\not=0$ for all $u\in\mathscr{H}$
since $K\not\equiv 0$. We set
\be\label{F}\mathcal{F}(u)=-\nabla \mathcal{J}(u)+\langle\nabla \mathcal{J}(u),\nabla \mathcal{L}(u)\rangle_{\mathscr{H}}\f{\nabla \mathcal{L}(u)}{\|\nabla \mathcal{L}(u)\|_{\mathscr{H}}^2}.\ee
By Lemma \ref{analytic} and the fact $\nabla\mathcal{L}(u)\not=0$ for all $u\in\mathscr{H}$,
we conclude that $\mathcal{F}\in C^\infty(\mathscr{H},\mathscr{H})$. Thus from the classical Cauchy-Lipschitz theorem
(\cite{Chang},  Theorem 1.9 of Chapter 1), there exists some $T>0$ such that
$u\in C^\infty([0,T);\mathscr{H})$ is a solution of
\be\label{equat}\le\{
\begin{array}{lll}
\p_tu=\mathcal{F}(u)\\[1.2ex]
u(0)=u_0\in\mathscr{S},
\end{array}
\ri.\ee
or equivalently (\ref{flow}).
In view of (\ref{F}), we have
$$\|\mathcal{F}(u)\|_{\mathscr{H}}\leq 2\|\nabla\mathcal{J}(u)\|_{\mathscr{H}}.$$
This together with (\ref{gradient-J}) leads to
$$\label{F-bd}\|\mathcal{F}(u)\|_{\mathscr{H}}\leq C\|u\|_{\mathscr{H}}+C.$$
Here and in the sequel, we often denote various constants by the same $C$. This together with the equation (\ref{equat}) implies that
$$\p_t\|u\|_{\mathscr{H}}^2=\langle\p_tu,u\rangle_{\mathscr{H}}\leq C\|u\|_{\mathscr{H}}^2+C,$$
which leads to
$$\p_t\le(e^{-Ct}\|u(t)\|_{\mathscr{H}}^2\ri)\leq Ce^{-Ct}.$$
Integrating this inequality from $0$ to $t<T$, one has
\be\label{bdd}\|u(t)\|_{\mathscr{H}}\leq (1+\|u_0\|_{\mathscr{H}})e^{C{T}/{2}}.\ee
It follows from (\ref{bdd}) that $u$ can be extended for all $t\in[0,\infty)$.

By (\ref{F}) and (\ref{equat}), we calculate
$$\p_t\mathcal{L}(u(t))=2\langle\nabla\mathcal{L}(u(t)),\p_tu\rangle_{\mathscr{H}}=
\langle\nabla\mathcal{L}(u(t)),\mathcal{F}(u)\rangle_{\mathscr{H}}=0.$$
Then we have for all $t\in[0,\infty)$,
$$\mathcal{L}(u(t))\equiv\mathcal{L}(u_0)=2\pi\chi(\Sigma,\bm\beta)$$
and thus $u(t)\in\mathscr{S}$. We now prove the energy identity (\ref{energy-id}). By (\ref{nab-C}),
$$\|\p_tu\|_{\mathscr{H}}^2=-\langle\nabla\mathcal{J}(u),\p_tu\rangle_{\mathscr{H}}
+\langle\nabla\mathcal{J}(u),\mathcal{N}(u)\rangle_{\mathscr{H}}\langle\mathcal{N}(u),\p_tu\rangle_{\mathscr{H}}.$$
Noting that
$$\langle\mathcal{N}(u),\p_tu\rangle_{\mathscr{H}}=\|\nabla\mathcal{L}(u)\|_{\mathscr{H}}^{-1}{\p_t\mathcal{L}(u)}=0,$$
we have
\be\label{der}\|\p_tu\|_{\mathscr{H}}^2=-\langle\nabla\mathcal{J}(u),\p_tu\rangle_{\mathscr{H}}=-\p_t\mathcal{J}(u).\ee
Integrating (\ref{der}) from $0$ to $t$, we obtain
$$\int_0^t\|\p_su(s)\|_{\mathscr{H}}^2ds=\mathcal{J}(u_0)-\mathcal{J}(u(t)).$$
This ends the proof of the Theorem.
$\hfill\Box$

\section{A sufficient condition for convergence}

In this section, we shall prove that if the solution $u(t)$ of (\ref{flow}) is uniformly bounded in $\mathscr{H}$, then
the flow must converge in $\mathscr{H}$. Precisely we have the following:

\begin{proposition}\label{prop1}
Let $u:[0,\infty)\ra\mathscr{H}$ be the solution of (\ref{flow}). Suppose that for all $t\in[0,\infty)$, there exists a
constant $C_0$ satisfying
\be\label{uniform}\|u(t)\|_{\mathscr{H}}\leq C_0.\ee
Then there exists some function $u_\infty\in W^{2,r}(\Sigma,g)\cap C^{\alpha}(\Sigma)$ for some $r>1$ and  $0<\alpha<1$, such that
$u(t)$ converges to $u_\infty$ in $\mathscr{H}$ as $t\ra\infty$. Moreover, if $\chi(\Sigma,\bm\beta)\not=0$, then
$u_\infty$ is a solution of (\ref{Kz-eqn}); if $\chi(\Sigma,\bm\beta)=0$, then $u_\infty+c$ is a solution of (\ref{Kz-eqn})
for some constant $c$.
\end{proposition}

\noindent{\it Proof.} By (\ref{energy-id}) and (\ref{uniform}), there exists a constant $C$ depending only on $C_0$ and $\kappa$
such that
$$\int_0^\infty\|\p_su(s)\|^2_{\mathscr{H}}ds\leq \mathcal{J}(u_0)+C.$$
As a consequence, there is a sequence $t_j\ra\infty$ satisfying
$$\label{tj-0}\|\p_tu(t_j)\|_{\mathscr{H}}=\|\nabla^{\mathscr{S}}\mathcal{J}(u(t_j))\|_{\mathscr{H}}\ra 0$$
as $j\ra\infty$. Since $\|u(t_j)\|_{\mathscr{H}}\leq C_0$ for all $j$, there would be some $u_\infty\in\mathscr{H}$ such that
up to a subsequence,
\bea
\label{weak}u(t_j)\rightharpoonup u_\infty&{\rm weakly\,\,in} &\mathscr{H}\\[1.2ex]\label{strong}
u(t_j)\ra u_\infty &{\rm strongly\,\,in}&L^q(\Sigma,g),\,\,\forall q>1.
\eea
Moreover, the singular Trudinger-Moser inequality (\cite{Troyanov}, Theorem 6) implies that
for any $\gamma>0$, there exists some constant $C$ depending only on $\gamma$ and the conical metric $g$ such that
\be\label{Tr-Mo}\int_\Sigma  e^{\gamma u(t_j)}dv_g\leq C.\ee

{\bf Claim 1.} {\it There holds $u_\infty\in\mathscr{S}$.}\\

To see this, we have by the mean value theorem
$$\int_\Sigma K(e^{2u(t_j)}-e^{2u_\infty})dv_g=\int_\Sigma Ke^\xi(2u(t_j)-2u_\infty)dv_g,$$
where $\xi$ lies between $2u(t_j)$ and $2u_\infty$. Clearly
$e^\xi\leq e^{2u(t_j)}+e^{2u_\infty}$. Thus in view of (\ref{Tr-Mo}), we estimate
\bna
\le|\int_\Sigma K(e^{2u(t_j)}-e^{2u_\infty})dv_g\ri|&\leq&2\sup_\Sigma |K|\le(\int_\Sigma(e^{4u(t_j)}+e^{4u_\infty})dv_g\ri)^{1/2}
\le(\int_\Sigma(u(t_j)-u_\infty)^2dv_g\ri)^{1/2}\\
&\leq& C\le(\int_\Sigma(u(t_j)-u_\infty)^2dv_g\ri)^{1/2}.
\ena
This together with (\ref{strong}) and the fact that $u_j\in\mathscr{S}$ leads to
$$\int_\Sigma K e^{2u_\infty}dv_g=\lim_{j\ra\infty}\int_\Sigma Ke^{2u(t_j)} dv_g=2\pi\chi(\Sigma,\bm\beta).$$
Hence $u_\infty\in \mathscr{S}$ and thus  Claim 1 follows.\\

{\bf Claim 2.} {\it There holds $\nabla^\mathscr{S}\mathcal{J}(u_\infty)=0$ and $u(t_j)\ra u_\infty$ in $\mathscr{H}$
as $j\ra\infty$.}\\

In view of (\ref{nab-C}), one has
\be\label{nab-C-j}\nabla^\mathscr{S}\mathcal{J}(u(t))=\nabla\mathcal{J}(u(t))-\langle\nabla\mathcal{J}(u(t)),
\nabla\mathcal{L}(u(t))\rangle_{\mathscr{H}}\f{\nabla\mathcal{L}(u(t))}
{\|\nabla\mathcal{L}(u(t))\|^2_{\mathscr{H}}}.\ee
We first prove that $\nabla^\mathscr{S}\mathcal{J}(u(t_j))$ converges to
$\nabla^\mathscr{S}\mathcal{J}(u_\infty)$ weakly in $\mathscr{H}$ as $j\ra\infty$.
To see this, it suffices to prove that as $j\ra\infty$,
\bea\label{J-infty}&&\nabla\mathcal{J}(u(t_j))\rightharpoonup\nabla\mathcal{J}(u_\infty)\quad{\rm weakly\,\,in}\quad\mathscr{H},\\
[1.2ex]\label{L-infty}&&\nabla\mathcal{L}(u(t_j))\rightharpoonup\nabla\mathcal{L}(u_\infty)\quad{\rm weakly\,\,in}\quad\mathscr{H},\\
[1.2ex]\label{inner-infty}&&\langle\nabla\mathcal{J}(u(t_j)),\nabla\mathcal{L}(u(t_j))\rangle_{\mathscr{H}}\ra
\langle\nabla\mathcal{J}(u_\infty),\nabla\mathcal{L}(u_\infty)\rangle_{\mathscr{H}},\\[1.2ex]\label{norm-infty}
&&\|\nabla\mathcal{L}(u(t_j))\|_{\mathscr{H}}\ra \|\nabla\mathcal{L}(u_\infty)\|_{\mathscr{H}}.
\eea
In view of (\ref{gradient-J}), we have
\be\label{nb-J}\nabla\mathcal{J}(u(t))=u(t)-(\Delta_g+I)^{-1}(u(t)-\kappa).\ee
For any $\phi\in\mathscr{H}$, one calculates
\bna
\langle(\Delta_g+I)^{-1}(u(t_j)+\kappa),\phi\rangle_{\mathscr{H}}&=&\int_\Sigma\nabla_g\le((\Delta_g+I)^{-1}(u(t_j)+\kappa)\ri)
\nabla_g\phi dv_g\\&&+\int_\Sigma(\Delta_g+I)^{-1}(u(t_j)+\kappa)\phi dv_g\\
&=&\int_\Sigma(\Delta_g+I)\le((\Delta_g+I)^{-1}(u(t_j)+\kappa)\ri)\phi dv_g\\
&=&\int_\Sigma(u(t_j)+\kappa)\phi dv_g.
\ena
This together with (\ref{weak}), (\ref{strong}) and (\ref{nb-J}) leads to (\ref{J-infty}).

In view of (\ref{gradient-2}),
\be\label{nb-L}\nabla\mathcal{L}(u(t))=(\Delta_g+I)^{-1}(Ke^{2u(t)}).\ee
For any $\phi\in\mathscr{H}$, one has as $j\ra\infty$,
\bna
\langle(\Delta_g+I)^{-1}(Ke^{2u(t_j)}),\phi\rangle_{\mathscr{H}}=\int_\Sigma Ke^{2u(t_j)}\phi dv_g
\ra\int_\Sigma Ke^{2u_\infty}\phi dv_g=\langle(\Delta_g+I)^{-1}(Ke^{2u_\infty}),\phi\rangle_{\mathscr{H}}.
\ena
This together with (\ref{nb-L}) leads to (\ref{L-infty}).

Let $f_j=(\Delta_g+I)^{-1}(Ke^{2u(t_j)})$, or equivalently $(\Delta_g+I)f_j=Ke^{2u(t_j)}$. Then standard elliptic estimates
lead to that $f_j$ is bounded in $W^{2,r}(\Sigma,g)$ for some $r>1$ and thus pre-compact in $\mathscr{H}$.
Up to a subsequence one may assume
$(\Delta_g+I)^{-1}(Ke^{2u(t_j)})$ converges to $(\Delta_g+I)^{-1}(Ke^{2u_\infty})$ in $\mathscr{H}$.
Similarly as before, one calculates
\bna
\langle\nabla\mathcal{J}(u(t_j)),\nabla\mathcal{L}(u(t_j))&=&\int_\Sigma\nabla_g(\Delta_g+I)^{-1}u(t_j)
\nabla_g(\Delta_g+I)^{-1}(Ke^{2u(t_j)}) dv_g\\&&+\int_\Sigma(\Delta_g+I)^{-1}u(t_j)(\Delta_g+I)^{-1}(Ke^{2u(t_j)}) dv_g\\
&=&\int_\Sigma u(t_j)(\Delta_g+I)^{-1}\le(Ke^{2u(t_j)}\ri) dv_g\\
&\ra&\int_\Sigma u_\infty(\Delta_g+I)^{-1}(Ke^{2u_\infty})dv_g\\
&=&\langle\nabla\mathcal{J}(u_\infty),\nabla\mathcal{L}(u_\infty)\rangle_{\mathscr{H}}.
\ena
This is exactly (\ref{inner-infty}). As for (\ref{norm-infty}), one has a strong estimate
\bea\nonumber
\|\nabla\mathcal{L}(u(t_j))\|_{\mathscr{H}}^2&=&\int_\Sigma Ke^{2u(t_j)}(\Delta_g+I)^{-1}(Ke^{2u(t_j)}) dv_g\\
\nonumber&\ra&\int_\Sigma Ke^{2u_\infty}(\Delta_g+I)^{-1}(Ke^{2u_\infty}) dv_g\\
&=&\|\nabla\mathcal{L}(u_\infty)\|_{\mathscr{H}}^2.\label{Str-L}
\eea
Therefore we have proved (\ref{J-infty})-(\ref{norm-infty}), and thus $\nabla^{\mathscr{S}}J(u(t_j))$
converges to $\nabla^{\mathscr{S}}J(u_\infty)$ weakly in $\mathscr{H}$. As a consequence
$$\|\nabla^{\mathscr{S}}J(u_\infty)\|_{\mathscr{H}}^2=\lim_{j\ra\infty}\langle\nabla^{\mathscr{S}}\mathcal{J}(u(t_j)),
\nabla^{\mathscr{S}}J(u_\infty)\rangle_{\mathscr{H}}\leq \lim_{j\ra\infty}\|\nabla^{\mathscr{S}}J(u(t_j))\|_{\mathscr{H}}
\|\nabla^{\mathscr{S}}J(u_\infty)\|_{\mathscr{H}}=0.$$
This immediately leads to $\nabla^{\mathscr{S}}J(u(t_j))$ converges in $\mathscr{H}$ to $\nabla^{\mathscr{S}}J(u_\infty)=0$.
It follows from (\ref{Str-L}) that $\nabla\mathcal{L}(u(t_j))$ converges in $\mathscr{H}$ to $\nabla\mathcal{L}(u_\infty)$.
Therefore, in view of (\ref{nab-C-j}) and (\ref{nb-J}), we obtain $u_j$ converges in $\mathscr{H}$ to $u_\infty$. This
concludes Claim 2.\\

By (\ref{nab-C-j}), (\ref{nb-J}) and (\ref{nb-L}), the equation $\nabla^{\mathscr{S}}J(u_\infty)=0$ is equivalent to
\be\label{lim}\Delta_gu_\infty+\kappa=c_\infty Ke^{2u_\infty}\ee
for some constant $c_\infty$. By elliptic estimates, we conclude that $u_\infty\in W^{2,r}(\Sigma,g)\cap C^\alpha(\Sigma)$ for some
$r>1$ and $0<\alpha<1$. If $\chi(\Sigma,\bm\beta)\not=0$, then we have by integrating (\ref{lim}), the Gauss-Bonnet formula and Claim 1
$$2\pi\chi(\Sigma,\bm\beta)=\int_\Sigma\kappa dv_g=c_\infty\int_\Sigma Ke^{2u_\infty}dv_g=2\pi\chi(\Sigma,\bm\beta) c_\infty.$$
It follows that $c_\infty=1$ and $u_\infty$ is a distributional solution of (\ref{Kz-eqn}). If $\chi(\Sigma,\bm\beta)=0$, then $\kappa=0$.
Multiplying both sides of (\ref{lim}) by $e^{-u_\infty}$, we have
$$-\int_\Sigma e^{-u_\infty}|\nabla_gu_\infty|^2dv_g=c_\infty\int_\Sigma Kdv_g,$$
which together with (\ref{compati}) implies that $c_\infty>0$. Then $u_\infty+\log c_\infty$ is a distributional solution of (\ref{Kz-eqn}).

Repeating the same argument of (\cite{BFR}, Pages 25-27), one can derive a Lojasiewicz-Simon inequality and then use it to obtain
$$\lim_{t\ra\infty}\|u(t)-u_\infty\|_{\mathscr{H}}=0.$$
This completes the proof of the proposition. $\hfill\Box$

\section{Convergence of the flow}

In this section, we prove Theorem \ref{theorem2} by using Proposition \ref{prop1}. The key point is to prove
that $\|u(t)\|_{\mathscr{H}}\leq C$ for all $t\in[0,\infty)$ under appropriate conditions.
\subsection{The null case}

{\it Proof of Theorem \ref{theorem2} in the null case.}
Suppose $\chi(\Sigma,\bm\beta)=0$. Since $\kappa$ is a constant, it follows from the Gauss-Bonnet formula that
$\kappa=0$. In view of (\ref{gradient-J}), on calculates
$$\Delta_gu(t)=(\Delta_g+I)\nabla\mathcal{J}(u(t)).$$
Integration by parts gives
$$\label{f}\int_\Sigma \nabla\mathcal{J}(u(t))dv_g=0,$$
which leads to
\be\label{jf}\langle \nabla\mathcal{J}(u(t)),1\rangle_{\mathscr{H}}=0.\ee
In view of (\ref{gradient-2}), we have
$$Ke^{2u(t)}=(\Delta_g+I)\nabla\mathcal{L}(u(t)).$$
Since $u(t)\in\mathscr{S}$, we have by integrating by parts
$$\int_\Sigma\nabla\mathcal{L}(u(t))dv_g=\int_\Sigma Ke^{2u(t)}dv_g=0.$$
Hence
\be\label{til-f}\langle\nabla\mathcal{L}(u(t)),1\rangle_{\mathscr{H}}=0.\ee
It follows from (\ref{jf}) and (\ref{til-f}) that
$$\p_t\int_\Sigma u(t)dv_g=\int_\Sigma \p_tudv_g=\langle\p_tu,1\rangle_{\mathscr{H}}=0.$$
Then there exists a constant $C$ such that
$$\label{int}\int_\Sigma u(t)dv_g\equiv C.$$
Using the Poincare inequality, we obtain
\be\label{l2}\int_\Sigma u^2dv_g\leq C\int_\Sigma|\nabla_gu|^2dv_g+C.\ee
By (\ref{energy-id}), there holds $\mathcal{J}(u(t))\leq \mathcal{J}(u_0)$, or equivalently
\be\label{leq}\int_\Sigma |\nabla_gu|^2dv_g\leq \int_\Sigma |\nabla_gu_0|^2dv_g.\ee
Combining (\ref{l2}) and (\ref{leq}), we obtain
$$\|u(t)\|_{\mathscr{H}}\leq C$$
for some constant $C$. This together with Proposition \ref{prop1} completes the proof of the theorem in the case $\chi(\Sigma,g)=0$.

\subsection{The negative case}

We first have a Poincar\'e inequality on conical surfaces.
\begin{lemma}\label{Poincare}
For all $u\in \mathscr{H}$, there holds
$$\int_\Sigma u^2dv_g\leq \f{1}{\lambda_g(\Sigma)}\int_\Sigma|\nabla_gu|^2dv_g+\f{1}{{\rm Vol}_g(\Sigma)}\le(\int_\Sigma udv_g\ri)^2,$$
where \be\label{eigen}\lambda_g(\Sigma)=\inf_{u\in \mathscr{H},\,\int_\Sigma udv_g=0,\, u\not\equiv 0}\f{\int_\Sigma|\nabla_gu|^2dv_g}
{\int_\Sigma u^2dv_g}.\ee
\end{lemma}

\noindent{\it Proof.} Applying a direct method of variation to (\ref{eigen}), one finds a function $u_0\in\mathscr{H}$ satisfying $\int_\Sigma u_0^2dv_g=1$ and
$$\lambda_g(\Sigma)=\int_\Sigma|\nabla_gu_0|^2dv_g>0.$$
Denote $$\overline{u}=\f{1}{{\rm Vol}_g(\Sigma)}\int_\Sigma udv_g.$$
By the definition of $\lambda_g(\Sigma)$, we have for all $u\in\mathscr{H}$,
$$\int_\Sigma|u-\overline{u}|^2dv_g\leq \f{1}{\lambda_g(\Sigma)}\int_\Sigma|\nabla_gu|^2dv_g.$$
Noting that
$$\int_\Sigma 2\overline{u}(u-\overline{u})dv_g=2\overline{u}\int_\Sigma (u-\overline{u})dv_g=0,$$
we obtain
\bna
\int_\Sigma u^2dv_g&=&\int_\Sigma\le((u-\overline{u})^2+\overline{u}^2+2\overline{u}(u-\overline{u})\ri)dv_g\\
&=&\int_\Sigma(u-\overline{u})^2dv_g+\overline{u}^2{\rm Vol}_g(\Sigma)\\
&\leq& \f{1}{\lambda_g(\Sigma)}\int_\Sigma|\nabla_gu|^2dv_g+\f{1}{{\rm Vol}_g(\Sigma)}\le(\int_\Sigma udv_g\ri)^2.
\ena
This gives the desired result. $\hfill\Box$\\

Next we have the following singular Trudinger-Moser inequality.

\begin{lemma}\label{Tr-Mo1}
There exist two constants $C$ and $\beta$ depending only on $(\Sigma,g)$ such  that for all
$u\in\mathscr{H}$,
\be\label{T-M}\int_\Sigma e^{2u}dv_g\leq C\exp\le(\beta\int_\Sigma|\nabla_gu|^2dv_g+\f{2}{{\rm Vol}_g(\Sigma)}\int_\Sigma udv_g\ri).\ee
\end{lemma}

\noindent{\it Proof.} Note that $g$ is a conical metric. The inequality (\ref{T-M}) follows from that of Troyanov
(\cite{Troyanov}, Theorem 6)
(see also Zhu \cite{Zhu-CMS} for a critical version). $\hfill\Box$\\

We remark that (\ref{T-M}) is a weak version of Trudinger-Moser inequality. For related strong versions, we refer the reader to
recent works \cite{Adi-Yang,Li-Yang,YangJDE,YangJGA,Yang-ZhuJFA} and the references therein.\\

\noindent{\it Proof of Theorem \ref{theorem2} in the negative case.}
Having Lemmas \ref{Poincare} and \ref{Tr-Mo1} in hand, we can prove an analog of (\cite{BFR}, Lemma 2)
by using the same method, and then repeating the argument of
the proof of
(\cite{BFR}, Part $(ii)$ of Theorem 2), we conclude the theorem in the case $\chi(\Sigma,\bm\beta)<0$. $\hfill\Box$\\

{\bf Acknowledgements}. This work is supported by National Science Foundation of China (Grant Nos. 11171347,
 11471014).

\end{document}